\documentclass[11pt,a4paper]{amsart} 
\usepackage[all]{xy}
\SelectTips{cm}{}
\calclayout \makeatletter
\def\serieslogo@{} 
\def\@setcopyright{} 
\makeatother
\title{A note on cohomological localization} 
\thanks{Version from October 2006.}
\dedicatory{Dedicated to the memory of Professor Paul M.~Cohn.}
\author[Henning Krause]{Henning Krause}
\address{Henning Krause\\ Institut f\"ur Mathematik\\
Universit\"at Paderborn\\ 33095 Paderborn\\ Germany.}
\email{hkrause@math.upb.de}

\theoremstyle{plain}
\newtheorem{lem}{Lemma}[section]
\newtheorem{prop}[lem]{Proposition}

\newtheorem{thm}[lem]{Theorem}

\theoremstyle{remark}
\newtheorem{rem}[lem]{Remark}

\theoremstyle{definition}

\newtheorem{defn}[lem]{Definition}

\numberwithin{equation}{section}

\newcommand{\id}{\operatorname{id}\nolimits}
\newcommand{\Id}{\operatorname{Id}\nolimits}

\newcommand{\End}{\operatorname{End}\nolimits}
\newcommand{\Hom}{\operatorname{Hom}\nolimits}

\newcommand{\Ker}{\operatorname{Ker}\nolimits}

\newcommand{\Ab}{\mathrm{Ab}}
\newcommand{\op}{\mathrm{op}}
\newcommand{\can}{\mathrm{can}}

\newcommand{\comp}{\mathop{\raisebox{+.3ex}{\hbox{$\scriptstyle\circ$}}}}
\newcommand{\lto}{\longrightarrow}
\newcommand{\xto}{\xrightarrow}

\def\a{\alpha}
\def\b{\beta}

\def\p{\phi}

\def\s{\sigma}

\def\Ga{\Gamma}
\def\La{\Lambda}
\def\Si{\Sigma}

\def\A{{\mathcal A}}

\def\C{{\mathcal C}}
\def\D{{\mathcal D}}

\def\X{{\mathcal X}}

\def\T{{\mathcal T}}
\def\U{{\mathcal U}}
\def\V{{\mathcal V}}

\def\bbZ{{\mathbb Z}}

\begin{document}

\begin{abstract}
Given a cohomological functor from a triangulated category to an
abelian category, we construct under appropriate assumptions for any
localization functor of the abelian category a lift to a localization
functor of the triangulated category. This discussion of Bousfield
localization is combined with a basic introduction to the concept of
localization for arbitrary categories.
\end{abstract}

\maketitle

\section{Introduction}
Localization is a technique to make morphisms in a category
invertible.  In this note we provide a basic introduction to this
concept and study a particular type of Bousfield localization. More
specifically, given a triangulated category $\T$ and a cohomological
functor $H^*\colon\T\to\A$ into an abelian category $\A$, we construct
(under appropriate assumptions) for any localization functor
$L\colon\A\to\A$ a new localization functor $\tilde L\colon \T\to\T$
such that $H^*\tilde L\cong LH^*$.

This paper is divided into two sections.  Localization functors are
introduced in full generality in the first section. This material is
classical. In fact, we present different approaches and explain how
they are related. Typical examples arise from localizations of
abelian and triangulated categories.  The second section is devoted to
studying a particular type of Bousfield localization. Here the main
result seems to be new.

The material of this note is based on some joint work with Dave Benson
and Srikanth Iyengar \cite{BIK}, where a special case of
Theorem~\ref{th:main} is used. Having in mind further applications, we
clarify results from this work by putting them into a broader
context. References to some of the classical sources, where localizations are
introduced and used, are listed at the end of this
note.

\section{Localization functors}
\subsection{Categories of fractions}
The following definition is taken from \cite{GZ}.
A functor $F\colon\C\to\D$ is said to make a morphism $\s$ of $\C$
invertible if $F\s$ is invertible. For a category $\C$ and any class
$\Si$ of morphisms of $\C$, there exists (ignoring set-theoretic
issues) the {\em category of fractions} $\C[\Si^{-1}]$ together with a
canonical functor
$$Q_\Si\colon\C\lto\C[\Si^{-1}]$$ having the following properties:
\begin{enumerate}
\item[(Q1)] $Q_\Si$ makes the morphisms in $\Si$ invertible.
\item[(Q2)] If a functor $F\colon\C\to\X$ makes the morphisms in $\Si$
invertible, then there is a unique functor $\bar F\colon \C[\Si^{-1}]\to\X$
such that $F=\bar F\comp Q_\Si$.
\end{enumerate}
Sometimes a functor of the form $Q_\Si\colon\C\to\C[\Si^{-1}]$ is
called localization functor. However, in this note we follow a
slightly different approach.  To be precise, we require the existence
of a right adjoint $G\colon\C[\Si^{-1}]\to\C$ and then we call the
composition $G\comp Q_\Si\colon\C\to\C$ a localization functor; see
Definition~\ref{de:loc}.

\subsection{Adjoint functors}

Let $F\colon\C\to\D$ and $G\colon\D\to\C$ be a pair of
functors and assume that $F$ is left adjoint to $G$. We denote by 
$$\Phi\colon F\comp G\to \Id\D\quad \text{and}\quad\Psi\colon\Id\C\to
G\comp F$$ the corresponding adjunction morphisms. Let $\Si$ denote
the class of morphisms $\s$  of $\C$ such that $F\s$ is invertible.

\begin{lem}\label{pr:quot}
The following statements are equivalent.
\begin{enumerate}
\item The functor $G$ is fully faithful.
\item The morphism $\Phi\colon F\comp G\to \Id\D$ is invertible.
\item The functor $\bar F\colon\C[\Si^{-1}]\to\D$ satisfying $F=\bar
F\comp Q_\Si$ is an equivalence.
\end{enumerate}
\end{lem}
\begin{proof}
See \cite[Sect.~I.1.3]{GZ}.
\end{proof}

\subsection{Localization functors}
Let $(F,G)$ be a pair of adjoint functors satisfying the equivalent
conditions of Lemma~\ref{pr:quot}. We consider the special
properties of the composition $L=G\comp F$, and it turns out that the pair
$(F,G)$ can be reconstructed from $L$ and the adjunction morphism
$\Psi\colon\Id\to L$.

\begin{lem}\label{le:loc}
Let $L\colon\C\to\C$ be a functor and $\Psi\colon\Id\C\to L$ be a morphism.
Then the following are equivalent.
\begin{enumerate}
\item $L\Psi\colon L\to L^2$ is invertible and $L\Psi=\Psi L$.
\item There exists an adjoint pair of functors $F\colon\C\to\D$ and
$G\colon\D\to\C$ such that $G$ is fully faithful, $L=G\comp F$, and
$\Psi\colon\Id\C\to G\comp F$ is the adjunction morphism.
\end{enumerate}
\end{lem}
\begin{proof}
(1) $\Rightarrow$ (2): Let $\D$ denote the full subcategory of $\C$
    formed by all objects $X$ such that $\Psi X$ is invertible. For each
    $X\in\D$, let $\Phi X\colon LX\to X$ be the inverse of $\Psi
    X$. Define $F\colon \C\to\D$ by $FX=LX$ and let $G\colon \D\to\C$
    be the inclusion. We claim that $F$ and $G$ form an adjoint pair. In
    fact, it is straightforward to check that the maps
$$\Hom_\D(FX,Y)\lto \Hom_\C(X,GY),\quad \a\mapsto G\a\comp\Psi X$$
and $$\Hom_\C(X,GY)\lto\Hom_\D(FX,Y),\quad \b\mapsto \Phi Y\comp F\b$$
are mutually inverse bijections.

(2) $\Rightarrow$ (1): Let $\Phi\colon FG\to\Id\D$ denote the second
adjunction morphism. Then the compositions $$F\xto{F\Psi}FGF\xto{\Phi
  F} F\quad\text{and}\quad G\xto{\Psi G}GFG\xto{G\Phi}G$$ are identity
morphisms. We know from Lemma~\ref{pr:quot} that $\Phi$ is invertible
because $G$ is fully faithful.  Therefore $L\Psi=GF\Psi$ is
invertible. Moreover, we have
$$L\Psi=GF\Psi=(G\Phi F)^{-1}=\Psi GF=\Psi L.$$
\end{proof}

\begin{rem}
Let $(L,\Psi)$ be a pair as above. Then the morphism $\Psi$ is
determined by $L$.  To see this, let $\D$ denote the full subcategory
of $\C$ formed by all objects isomorphic to one of the form $LX$ with
$X$ in $\C$. Then the inclusion $G\colon \D\to\C$ has a left adjoint
$F$ taking $X\in\C$ to $LX$, and $\Psi$ equals the adjunction
morphism $\Id\C\to GF=L$.
\end{rem}

\begin{defn}
\label{de:loc}
A functor $L\colon\C\to\C$ is called \emph{localization functor} if
there exists a morphism $\Psi\colon\Id\C\to L$ such that the pair
$(L,\Psi)$ satisfies the equivalent conditions of Lemma~\ref{le:loc}.
We call $\Psi$ the \emph{adjunction morphism} for $L$.
\end{defn}

\subsection{Local objects}
Let $L\colon\C\to\C$ be a localization functor and $\Psi\colon\Id\C\to
L$ its adjunction morphism.  We write $\Si_L$ for the class of
morphisms of $\C$ which $L$ makes invertible.

\begin{lem} 
\label{le:locobj}
The following statements are equivalent for an object $X$ in $\C$.
\begin{enumerate}
\item $\Psi X\colon X\to LX$ is invertibe.
\item $X\cong LX'$ for some $X'\in\C$.
\item For all $\s\colon V\to W$ in $\Si_L$, the induced map
$$\Hom_\T(\s,X)\colon\Hom_\C(W,X)\to\Hom_\T(V,X)$$ is bijective.
\end{enumerate}
\end{lem}
\begin{proof}
(1) $\Rightarrow$ (2): Take $X'=X$.

(2) $\Rightarrow$ (3): Fix $\s\colon V\to W$ and consider the
    following commutative diagram.
$$\xymatrix{\Hom_\T(LW,LX)\ar[d]^-{(L\s,LX)}\ar[rr]^-{(\Psi
W,LX)}&&\Hom_\T(W,LX)\ar[d]^-{(\s,LX)}\\
\Hom_\T(LV,LX)\ar[rr]^-{(\Psi V,LX)}&&\Hom_\T(V,LX) }$$
The horizontal maps are invertible. It follows that 
$\Hom_\T(\s,LX)$ is invertible, provided that $L\s$ is invertible.

(3) $\Rightarrow$ (1): The morphism $\Psi X$ belongs to $\Si_L$
because $L$ is a localization functor. Thus the assumption implies
that the induced map $\Hom_\C(LX,X)\to\Hom_\T(X,X)$ is
bijective. Let $\a\colon LX\to X$ be the morphism which is sent to $\id X$. Then
$\a$ is an inverse of $\Psi X$, and therefore $\Psi X$ is invertible.
\end{proof}
\begin{defn} 
An object $X$ in $\C$ is called {\em $L$-local} if it satisfies the
equivalent conditions of Lemma~\ref{le:locobj}.
\end{defn}

We denote by $\C_L$ the full subcategory of $\C$ which is formed by
all $L$-local objects. Observe that $\C_L$ determines $L$ up to an
isomorphism. This is a consequence of the first part of the following
proposition, because a left adjoint of the inclusion $\C_L\to\C$ is
unique up to an isomorphism.

\begin{prop}
\label{pr:loc}
Let $L\colon\C\to\C$ be a localization functor.
\begin{enumerate}
\item The functor $\C\to \C_L$ taking $X\in\C$ to $LX$ is a left
adjoint for the inclusion functor $\C_L\to\C$. 
\item The functor $L$ induces an equivalence
$$\C[\Si_L^{-1}]\stackrel{\sim}\lto\C_L.$$
\end{enumerate}
\end{prop}
\begin{proof}
The first part follows from Lemma~\ref{le:loc} and its proof. The
second part then follows from Lemma~\ref{pr:quot}.
\end{proof}

\subsection{Acyclic objects}
Let $\C$ be an additive category and $L\colon\C\to\C$ be an additive
localization functor.  An object $X$ in $\C$ is called {\em $L$-acyclic} if
$LX=0$. We denote by $\Ker L$ the full subcategory formed by all
$L$-acyclic objects.

\begin{lem} \label{le:acyclic}
An object $X$ in $\C$ is $L$-acyclic if and only if $\Hom_\C(X,Y)=0$
for all $L$-local objects $Y$.
\end{lem}
\begin{proof}
Write $L=G\comp F$ as the composite of two adjoint
functors; see Lemma~\ref{le:loc}. We have for each pair $X,Y$ in $\C$
$$\Hom_\T(LX,LY)\cong\Hom_\T(FX,FY)\cong\Hom_\T(X,LY).$$ Now the
assertion follows from the fact $Y\cong LY$ for each $L$-local object
$Y$.
\end{proof}

\subsection{Localization for abelian categories}
Let $\A$ be an abelian category and let $L\colon\A\to\A$ be an exact
localization functor. Then $\A_L$ and $\Ker L$ are abelian categories.
Moreover, the inclusion $\Ker L\to\A$ and the functor $\A\to \A_L$
sending $X\in\A$ to $LX$ are both exact. Note that $L$ induces an equivalence
$$\A/{\Ker L}\stackrel{\sim}\lto\A_L$$ where $\A/{\Ker L}$ denotes the
quotient category in the sense of Gabriel \cite{G}.  Recall that
$\A/{\Ker L}=\A[\Si^{-1}]$, where $\Si$ denotes the class of morphisms of
$\A$ such that its kernel and cokernel belong to $\Ker L$.  The
equivalence $\A/{\Ker L}\stackrel{\sim}\to\A_L$ now follows from
Proposition~\ref{pr:loc} since $\Si$ equals the class of morphisms of $\A$
which $L$ makes invertible.

\subsection{Localization for triangulated categories}
Let $\T$ be a triangulated category and let $L\colon\T\to\T$ be an
exact localization functor. Then $\T_L$ and $\Ker L$ are triangulated
categories.  Moreover, the inclusion $\Ker L\to\T$ and the functor
$\T\to \T_L$ sending $X\in\T$ to $LX$ are both exact. Note that $L$ induces
an equivalence
$$\T/{\Ker L}\stackrel{\sim}\lto\T_L$$
where $\T/{\Ker L}$ denotes the quotient category in the sense of Verdier \cite{V}.
Recall that
$\T/{\Ker L}=\T[\Si^{-1}]$, where $\Si$ denotes the class of morphisms of
$\T$ such that its cone belongs to $\Ker L$.  The
equivalence $\T/{\Ker L}\stackrel{\sim}\to\T_L$ now follows from
Proposition~\ref{pr:loc} since $\Si$ equals the class of morphisms of $\T$
which $L$ makes invertible.

\section{Cohomological localization}
Let $\T$ be a triangulated category which admits small coproducts. We
assume that $\T$ is generated by a set of compact objects. Let $S$
denote the suspension functor of $\T$.

We fix a graded\footnote{All graded rings and modules are graded over
$\bbZ$. Morphisms between graded modules are degree zero maps.} ring
$\La$ and a graded cohomological functor
$$H^*\colon\T\lto\A$$ into the category $\A$ of graded $\La$-modules.
Thus $H^*$ is a functor which sends each exact triangle in $\T$ to an
exact sequence in $\A$, and we have an isomorphism $H^*\comp S\cong
T\comp H^*$ where $T$ denotes the shift functor for $\A$.  In addition, we
assume that $H^*$ preserves small products and coproducts.

\begin{thm}
\label{th:main}
Let $(L,\Psi)$ be an exact localization functor for the category $\A$
of graded $\La$-modules. Then there exists an exact localization functor
$(\tilde L,\tilde\Psi)$ for $\T$ such that the following diagram
commutes up to a natural isomorphism.
$$\xymatrix{ \T\ar[rr]^-{\tilde L}\ar[d]^-{H^*}&&\T\ar[d]^-{H^*}\\
\A\ar[rr]^-L&&\A }$$ More precisely, we have the following.
\begin{enumerate}
\item The morphisms $LH^*\tilde\Psi$, $\Psi H^*\tilde L$, and
$$LH^*\xto{LH^*\tilde\Psi}LH^*\tilde L\xto{(\Psi H^*\tilde L)^{-1}}H^*\tilde L$$ are
invertible.
\item An object $X$ in $\T$ is $\tilde L$-acyclic if and only if $H^*X$
is $L$-acyclic.
\item If an object $X$ in $\T$ is $\tilde L$-local, then $H^*X$ is
$L$-local. The converse holds, provided that $H^*$ reflects
isomorphisms.
\end{enumerate}
\end{thm}

A special case of this theorem with almost identical proof appears in
\cite{BIK}.

\subsection{Cohomological functors}
Recall that a functor from $\T$ to any abelian category is
\emph{cohomological} if it sends each exact triangle in $\T$ to an
exact sequence.  Our assumption on $\T$ implies that $\T$ satisfies
Brown representability: each cohomological functor $\T^\op\to\Ab$
sending small coproducts in $\T$ to products is representable, and
each cohomological functor $\T\to\Ab$ sending small products in $\T$
to products is representable; see \cite{N}.

We consider some cohomological functors in more detail.
Given objects $X,Y$ in $\T$, we obtain a graded abelian group
\[ \Hom^*_\T(X,Y)=\coprod_{i\in\bbZ}\Hom_\T(X,S^i Y).\]
Note that $\End^*_\T(X)=\Hom^*_\T(X,X)$ has a natural structure of a
graded ring and $\Hom^*_\T(X,Y)$ is a graded $\End^*_\T(X)$-module.

Now fix an object $C$ in $\T$ and a homomorphism of graded rings
$\Ga\to \End_\T^*(C)$.  For each $X$ in $\T$, we view $\Hom^*_\T(C,X)$
as a graded $\Ga$-module and obtain a graded cohomological functor
from $\T$ into the category of graded $\Ga$-modules. It turns out that
the functor $H^*$ is precisely of this form.

\begin{lem}
\label{le:fun}
There exists a compact object $C$ in $\T$, a graded ring homomorphism
$\La\to\End_\T^*(C)$, and a natural isomorphism
$H^*X\cong\Hom^*_\T(C,X)$ of graded $\La$-modules for all $X$ in $\T$.
\end{lem}
\begin{proof}
Let $H^nX$ denote the degree $n$ part of $H^*X$. Using Brown
representability, we find a compact object $C$ in $\T$ representing
the functor $H^0\colon\T\to\Ab$. Thus $H^0X\cong\Hom_\T(C,X)$ for all
$X$ in $\T$, and this extends to an isomorphism
$H^*X\cong\Hom^*_\T(C,X)$ since $H^*(S^nX)\cong (H^*X)[n]$ for all
$n\in\bbZ$.  We obtain a graded ring homomorphism $\La\to\End_\T^*(C)$
because $H^*$ is a functor into graded $\La$-modules. Thus
$H^*X\cong\Hom^*_\T(C,X)$ becomes an isomorphism of $\La$-modules.
\end{proof}

\begin{rem}
The functor $H^*\cong\Hom_\T^*(C,-)$ reflects isomorphism if and only if
$C$ generates $\T$. Recall that a compact object $C$ \emph{generates} $\T$
if $\Hom_\T^*(C,X)=0$ implies $X=0$ for all objects $X$ in $\T$.
\end{rem}

\subsection{Proof of the theorem}
We use Lemma~\ref{le:fun} and identify $H^*=\Hom^*_\T(C,-)$ for some
appropriate object $C$ in $\T$.

\begin{proof}[Proof of Theorem~\ref{th:main}]
Let $\A_L$ denote the full subcategory of $\A$ formed by all $L$-local
objects. This subcategory is {\em coherent} (that is, for any exact
sequence $X_1\to X_2\to X_3\to X_4\to X_5$ with
$X_1,X_2,X_4,X_5\in\A$, we have $X_3\in\A$), and closed under taking
small products. The $L$-local objects form a Grothendieck category and
therefore $\A_L$ admits an injective cogenerator, say $I$; see
\cite{G}.  Using Brown's representability theorem, there exists
$\tilde I\in\T$ such that
\begin{equation}\label{eq:1}\Hom_\A(H^*-,I)\cong\Hom_\T(-,\tilde I).
\end{equation} Note that this extends to
an isomorphism
$$\Hom^*_\A(H^*-,I)\cong\Hom^*_\T(-,\tilde I).$$ Now consider the
subcategory $\V$ of $\T$ which is formed by all objects $X\in\T$ such
that $H^*X$ is $L$-local. This is a triangulated subcategory
which is closed under taking small products.  Observe that $\tilde I$
belongs to $\V$. To prove this, take a free presentation
$$F_1\lto F_0\lto H^*C\lto 0$$ over $\La$ and apply $\Hom_\A^*(-,I)$ to
it. Using the isomorphism \eqref{eq:1}, we see that
$H^*\tilde I$ belongs to $\A_L$ because $\A_L$ is coherent and closed
under taking small products.

Now let $\U$ denote the smallest triangulated subcategory of $\T$
containing $\tilde I$ and closed under taking small products. Observe
that $\U\subseteq \V$. We claim that $\U$ is perfectly cogenerated by
$\tilde I$ in the sense of \cite{K}. Thus, given a family of morphisms
$\p_i\colon X_i\to Y_i$ in $\U$ such that $\Hom_\T(Y_i,\tilde
I)\to\Hom_\T(X_i,\tilde I)$ is surjective for all $i$, we need to show
that $\Hom_\T(\prod_iY_i,\tilde I)\to\Hom_\T(\prod_iX_i,\tilde I)$ is
surjective. We argue as follows. If $\Hom_\T(Y_i,\tilde
I)\to\Hom_\T(X_i,\tilde I)$ is surjective, then the isomorphism
\eqref{eq:1} implies that $H^*\p_i$ is a monomorphism since $I$ is an
injective cogenerator for $\A_L$. Thus the product $\prod_i\p_i\colon
\prod_iX_i\to \prod_iY_i$ induces a monomorphism
$H^*\prod_i\p_i=\prod_iH^*\p_i$ and therefore
$\Hom_\T(\prod_i\p_i,\tilde I)$ is surjective.  We conclude from
Brown's representability theorem \cite{K} that the inclusion functor
$G\colon \U\to\T$ has a left adjoint $F\colon\T\to\U$. Now put $\tilde
L=G\comp F$ and let $\tilde\Psi\colon\Id\T\to\tilde L$ denote the
adjunction morphism.  Then $(\tilde L,\tilde\Psi)$ is a localization
functor by Lemma~\ref{le:loc}.

Next we show that an object $X\in\T$ is $\tilde L$-acyclic if and only
if $H^*X$ is $L$-acyclic. This follows from Lemma~\ref{le:acyclic}
and the isomorphism \eqref{eq:1}, because we have
$$\tilde LX=0 \iff \Hom_\T(X,\tilde I)=0\iff\Hom_\A(H^*X, I)=0\iff LH^*X=0.$$

Now consider the following commutative square.
\begin{equation}\label{eq:2}
\xymatrix{H^*X\ar[rr]^-{\Psi H^*X}\ar[d]^-{H^*\tilde\Psi X}
&&LH^*X\ar[d]^-{LH^*\tilde\Psi X}\\ H^*\tilde LX\ar[rr]^-{\Psi
H^*\tilde LX}&&LH^*\tilde LX}
\end{equation} 
We need to show that $LH^*\tilde\Psi X$ and $\Psi H^*\tilde L X$ are
invertible for each $X$ in $\T$. The morphism $\tilde\Psi X$ induces an
exact triangle
$$X'\to X\xto{\tilde\Psi X}\tilde LX\to S X'$$ with $\tilde
LX'=0=\tilde LS X'$. Applying the cohomological functor $LH^*$, we
see that $LH^*\tilde\Psi X$ is an isomorphism, since $LH^*X'=0=LH^*S
X'$.  Thus $LH^*\tilde\Psi$ is invertible.  The morphism $\Psi H^*\tilde
LX$ is invertible because $H^*\tilde LX$ is $L$-local. This follows
from the fact that $\tilde LX$ belongs to $\U$.

The commutative square \eqref{eq:2} implies that $H^*\tilde\Psi X$ is
invertible if and only if $\Psi H^*X$ is invertible.  Thus if $X$ is
$\tilde L$-local, then $H^*X$ is $L$-local. The converse holds if
$H^*$ reflects isomorphisms. 
\end{proof}

\begin{rem}
The localization functor $\tilde L$ is unique up to isomorphism. More
precisely, $\tilde L$ is determined by $H^*$ and $L$ via part (2) of the theorem,
since any localization functor is determined by its class of acyclic
objects.
\end{rem}

\begin{rem}
The composite $L\comp H^*\colon\T\to\A$ induces a functor
$\T/{\Ker\tilde L}\to\A/{\Ker L}$ making the following diagram
commutative.
$$\xymatrix{ \T\ar[rr]^-{\can}\ar[d]^{H^*}&&\T/{\Ker\tilde
L}\ar[rr]^-\sim\ar[d]^{J}&&\T_{\tilde L}\ar[d]^{H^*}\\
\A\ar[rr]^-{\can}&&\A/{\Ker L}\ar[rr]^-\sim&&\A_L }$$ The functor $J$
has an explicit description. For this, we may assume that the
canonical functors induce identity maps on objects. Given an object $X$
in $\T$, we have the following chain of natural isomorphisms.
$$\Hom^*_{\T/{\Ker\tilde L}}(C,X)\cong\Hom_{\T_{\tilde
L}}^*(\tilde LC,\tilde LX)\cong\Hom_\T^*(C,\tilde LX)=H^*\tilde LX\cong LH^*X$$
Thus we can identify $$JX=\Hom^*_{\T/{\Ker\tilde L}}(C,X)$$ where the graded
abelian group $\Hom^*_{\T/{\Ker\tilde L}}(C,X)$ is viewed as a graded
$\La$-module via the action of $\End^*_{\T/{\Ker\tilde L}}(C)$ and the
homomorphism
$$\La\lto\End^*_{\T}(C)\lto\End^*_{\T/{\Ker\tilde L}}(C).$$
\end{rem}

\begin{rem}
Suppose that $C$ is a generator of $\T$. If $L$ preserves small
coproducts, then it follows that $\tilde L$ preserves small
coproducts.  In fact, the assumption implies that $H^*\tilde L$
preserves small coproducts, since $LH^*\cong H^*\tilde L$. But $H^*$
reflects isomorphisms because $C$ is a generator of $\T$. Thus $\tilde
L$ preserves coproducts.
\end{rem}

\end{document}